\documentclass[review]{elsarticle}

\usepackage{lineno,hyperref}
\modulolinenumbers[5]

\usepackage{amssymb, amsfonts, amsthm, amsmath, cite}
\usepackage{graphicx}

\journal{Journal of \LaTeX\ Templates}

\begin{document}

\begin{frontmatter}

\title{A NEW METHOD IN THE PROBLEM \\OF THREE CUBES}

\author{Armen Avagyan, Gurgen Dallakyan}
\address{Armenian State Pedagogical University after Khachatur Abovyan}


\ead{avagyana73@gmail.com}


\begin{abstract}
In the current paper we are seeking $P_1(y), P_2(y), P_{3}(y)$  with the highest possible degree polynomials with integer coefficients, and $Q(y)$  via the lowest possible degree polynomial, such that  $P_1^3(y)+P_2^3(y)+P_3^3(y)=Q(y)$. Actually, the solution of this problem has close relation with the problem of the sum of three cubes  $a^3+b^3+c^3=d$, since $\deg Q(y)=0$  case coincides with above mentioned problem. It has been considered estimation of possibility of minimization of  $\deg Q(y)$. 

As a conclusion, for specific values of $d$  we survey a new algorithm for finding integer solutions of  $a^3+b^3+c^3=d$.
\end{abstract}


\end{frontmatter}


\section{\large Some facts on the history of the equation $a^3+b^3+c^3=d$}	

The question as to which integers are expressible as a sum of three integer cubes is over 160 years old. The first known reference to this problem was made by Fermat, who has offered to find the three nonzero integers, so that the sum of their nth powers is equal to zero. This framework of the problem makes the beginning of survey of  equation $a^3+b^3+c^3=d$ for many mathematicians. We present some significant result schedule. 
1825 year -- S. Ryley in [1] gave a parametrization of rational solutions for $d\in Z$:
$$x=\frac{(9k^6-30d^2 k^3+d^4)(3k^3+d^2)+72d^4 k^3}{6kd (3k^3+d^2)^2}$$
$$y=\frac{30d^2 k^3-9k^6-d^4}{6kd (3k^3+d^2)}$$
$$z=\frac{18d k^5-6d^3 k^2}{(3k^3+d^2)^2}$$
1908 year -- A.S. Werebrusov [2] found the following parametric family for $d=2$:
$$(6t^3+1)^3-(6t^3-1)^3-(6t^2)^3=2$$
1936 year -- Later in [3] Mahler  discovered a first parametric solution for $d=1$:
$$(9t^4)^3+(3t-9t^4)^3+(1-9t^3)^3=1$$

\hangindent=0.8cm
\noindent 
1942 year -- Mordell proved in [2] that for any other $d$ a parametric solution with rational coefficients must have degree at least 5.

\hangindent=0.8cm
\noindent 
1954 year -- Miller and Woollett [4] discovered explicit representations for 69 values of  $d$ between 1 and 100. Their search exhausted the region $\vert a\vert, \vert b \vert, \vert c \vert \leqslant  3164$.

\hangindent=0.8cm
\noindent 
1963 year -- 1963 Gardiner, Lazarus, and Stein [5] looked at the equation $x^3+y^3=z^3-d$ in the range 
$0\leqslant x\leqslant y\leqslant 2^{16}$, where $0\leqslant z-x\leqslant 2^{16}$ and $0\leqslant \vert d \vert \leqslant 999$. Their search left only 70 values of $d$ between 1 and 1000 without a known representation including eight values less than 100.

\hangindent=0.8cm
\noindent 
1992 year -- the first solution for $d=39$ was found. Heath-Brown, Lioen, and te Riele [8] determined that $39=134476^3+117367^3+(-159380)^3$ with the rather deep algorithm of Heath-Brown [6]. This algorithm involved searching for solutions  for a specific value of $d$ using the class number of $Q(\sqrt[3]{d})$ to eliminate values of $a$, $b$, $c$ which would not yield a solution.

\hangindent=0.8cm
\noindent 
1994 year -- Koyama [7] used modern computers to expand the search region to $\vert a \vert , \vert b \vert , \vert c \vert\leqslant 2^{21}$ and successfully found first solutions for 16 integers between 100 and 1000 [9]. Also in 1994, Conn and Vaserstein [8] chose specific values of $d$ to target, and then used relations implied by each chosen value to limit the number of triples $(a, b, c)$ searched.  So doing, they found first representations for 84 and 960. Their paper also lists a solution for each $d<100$ for which a representation was known.

\hangindent=0.8cm
\noindent 
1995 year -- Bremner [9] devised an algorithm which uses elliptic curve arguments to narrow the search space. He discovered a solution for 75 (and thus a solution for 600), leaving only five values less than 100 for which no solution was known. Lukes then extended this search method to also find the first representations for each of the values 110, 435, and 478 [10]. 

\hangindent=0.8cm
\noindent 
1997 year -- Koyama, Tsuruoka, and Sekigawa [11] used a new algorithm to find first solutions for five more values between 100 and 1000 as well as independently finding the same solution for 75 that Bremner found. Also in the same paper, the authors discuss the complexity of the above algorithms.

\hangindent=0.8cm
\noindent 
1999 year -- Bernstein [12] had implemented the method of Elkies [13] and found solutions for 11 new values of $d$.

Summarizing the above, it can be noted that up to 21st century only 27 values were left unresolved. These, together with the range of their search, are presented in the table below.

\begin{center}
\begin{tabular}{|c|c|}
\hline
$k$ & $< T$\\
\hline
33 & $10^{\scriptscriptstyle{\wedge}}{12}$\\
\hline
42 & $6.5\times 10^{\scriptscriptstyle{\wedge}}{11}$\\
\hline
74 & $1.5\times 10^{\scriptscriptstyle{\wedge}}{11}$\\
\hline
156, 165, 318, 366, 390, 420, 534, 564, 579, & \\
609, 627, 633, 732, 758, 786, 789, 795, 834, & $10^{10}$\\
894, 903, 906, 921, 948, 975 & \\
\hline
\end{tabular}
\end{center}

Only recently, in 2007, Elsenhans and Yahnel [14] found the solutions for the values \\
$d=156,\,318,\,366,\,420,\,564,\,758,\,789,\,894,\,948$.\\
$156= 26577110807569^{\scriptscriptstyle{\wedge}}3 -18161093358005^{\scriptscriptstyle{\wedge}}3 -23381515025762^{\scriptscriptstyle{\wedge}}3$\\
$318= 47835963799^{\scriptscriptstyle{\wedge}}3+ 20549442727^{\scriptscriptstyle{\wedge}}3 -49068024704^{\scriptscriptstyle{\wedge}}3$\\
$318= 1970320861387^{\scriptscriptstyle{\wedge}}3+ 1750553226136^{\scriptscriptstyle{\wedge}}3 -2352152467181^{\scriptscriptstyle{\wedge}}3$\\
$318= 30828727881037^{\scriptscriptstyle{\wedge}}3+ 27378037791169^{\scriptscriptstyle{\wedge}}3 -36796384363814^{\scriptscriptstyle{\wedge}}3$\\
$366= 241832223257^{\scriptscriptstyle{\wedge}}3+ 167734571306^{\scriptscriptstyle{\wedge}}3 -266193616507^{\scriptscriptstyle{\wedge}}3$\\
$420= 8859060149051^{\scriptscriptstyle{\wedge}}3 -2680209928162^{\scriptscriptstyle{\wedge}}3 -8776520527687^{\scriptscriptstyle{\wedge}}3$\\
$564= 53872419107^{\scriptscriptstyle{\wedge}}3 -1300749634^{\scriptscriptstyle{\wedge}}3 -53872166335^{\scriptscriptstyle{\wedge}}3$\\
$758= 662325744409^{\scriptscriptstyle{\wedge}}3+ 109962567936^{\scriptscriptstyle{\wedge}}3 -663334553003^{\scriptscriptstyle{\wedge}}3$\\
$789= 18918117957926^{\scriptscriptstyle{\wedge}}3+ 4836228687485^{\scriptscriptstyle{\wedge}}3 -19022888796058^{\scriptscriptstyle{\wedge}}3$\\
$894= 19868127639556^{\scriptscriptstyle{\wedge}}3+ 2322626411251^{\scriptscriptstyle{\wedge}}3 -19878702430997^{\scriptscriptstyle{\wedge}}3$\\
$948= 323019573172^{\scriptscriptstyle{\wedge}}3+ 63657228055^{\scriptscriptstyle{\wedge}}3 -323841549995^{\scriptscriptstyle{\wedge}}3$\\
$948= 103458528103519^{\scriptscriptstyle{\wedge}}3+ 6604706697037^{\scriptscriptstyle{\wedge}}3 -103467499687004^{\scriptscriptstyle{\wedge}}3$
\\
Thus, until 1000 there are only the numbers 33, 42, 114, 165, 390, 579, 627, 633, 732, 795, 906, 921, 975 lefts, which have not yet been solved, all the other presentations are posted on the web, in particular, it was made by Sander Huisman [15].

\section{\large Some known notes about the equation $a^3+b^3+c^3=d$ for specific values of $d$}

Consider first the equation
\begin{equation}\label{2-1}
a^3+b^3+c^3=1.                                                 
\end{equation}
This has infinitely many solutions because of the identity
\begin{equation}\label{2-2}
(1+-9m^3)^3+(9m^4)^3+(-9m^4-+3m)^3=1
\end{equation}
but there are other solutions as well.  Are there any other identities that give a different 1-parameter family of solutions?  Is every solution of \eqref{2-1} a member of a family like this? In general it's known that there is no finite method for determining whether a given Diophantine equation has solutions.  However, it’s an open problem whether if there is a general method for determining if a given Diophantine equation has ''algebraic'' solutions, i.e., an algebraic identity like the one above that gives an infinite family of solutions.  More specifically, is there a proposition, that only equations of    $genus<2$ can have an algebraic solution.

It may be worth mentioning that the complete rational-solution of the equation $a^3+b^3+c^3=d^3$ is known, and is given by
\begin{align*}
& a=q[1-(x-3y)(x^2+3y^2)]\\
& b=-q[1-(x+3y)(x^2+3y^2)]\\
& c=q[(x^2+3y^2)^2-(x+3y)]\\
& d=q[(x^2+3y^2)^2-(x-3y)]
\end{align*}
where $q$, $x$, $y$ are any rational numbers. So if we set $q$ equal to the inverse of  
$[(x^2+3y^2)^2-(x-3y)]$ we have rational solutions of \eqref{2-1}.
  
However, the problem of finding the integer-solutions is more difficult.  If $d$ is allowed to be any integer (not just 1) then Ramanujan gave the integer solutions
\begin{align*}
& a=3{{n}^{2}}+5nm-5{{m}^{2}},\\
& b=4{{n}^{2}}-4nm+6{{m}^{2}},\\
& c=5{{n}^{2}}-5nm-3{{m}^{2}},\\
& d=6{{n}^{2}}-4nm+4{{m}^{2}}.
\end{align*}

This occasionally gives a solution of equation \eqref{2-1} (with appropriate changes in sign), as in the following cases
\begin{center}
\begin{tabular}{ccccc}
n & m  &   a & b & c \\
\hline
1   &  -1   &   1  &   2    &   -2\\
1  &   -2   &    9   &   10   &  -12\\
5  &  -12  &  -135 &  -138 &  172\\
19  &   -8  &  -791 &  -812 &  1010\\
46  & -109 &  11161 &  11468 &  -14258\\
73 &  -173 &  65601 &   67402&  -83802\\
419  & -993 &  -951690 & -926271 &  1183258\\
\end{tabular}
\end{center}
However, this doesn't cover all of the solutions given by \eqref{2-2}. By the way, the equation  
$a^3+b^3+c^3=1$ has algebraic solutions [16], other than \eqref{2-2}.

There are  known to be infinitely many algebraic solutions, for example
$$(1-9t^3+648t^6+3888t^9)^3+(-135t^4+3888t^{10})^3+(3t-81t^4-1296t^7-3888t^{10})^3=1$$
However, it's not known whether every solution of the equation lies in some family of solutions with an algebraic parameterization.

Interestingly, note that if you replace 1 by 2, then again there's a parametric solution:
\begin{equation}\label{2-3}
(6t^3+1)^3-(6t^3-1)^3-(6t^2)^3=2                               
\end{equation}                        
but again this doesn't cover all known integer solutions. Note, that precisely one solution is known that is not given by \eqref{2-3} (see~[16]):
$$1214928^3+3480205^3-3528875^3=2$$
It's evidently not known up todays if there are any other algebraic solutions besides the one noted above. 

In general it seems to be a difficult problem to characterize all the solutions of                         
\begin{equation}\label{2-4}
a^3+b^3+c^3=d
\end{equation}
for some arbitrary integer $d>2$.  In particular, the question of whether all integer solutions are given by an algebraic identity seems both difficult and interesting.

Note that for $d\equiv \pm 4(\bmod 9)$ there are no solutions since, for any integer $a$, 
$a^3\equiv 0, 1, -1(\bmod 9)$. 
It is a long standing problem as to whether every rational integer $d\neq 4, 5(\bmod 9)$ can be written as a sum of three integral cubes. According to the web page [12] of Daniel Bernstein, the first attacks by computer were carried out as early as 1955.

Nevertheless, for example, for $d=3$, there is still no solution known apart from the obvious ones: $(1, 1, 1)$, $(4, 4,-5)$, $(4,-5, 4)$, and $(-5, 4, 4)$. For $d=30$, the first solution was found by N.~Elkies and his coworkers in 2000 [17]. 
It is interesting to note that, in 1992, D.R.~Heath-Brown [6] had made a prediction on the density of the solutions for $d=30$ without knowing any solution explicitly.

Over the years, a number of algorithms have been developed in order to attack the general problem. An excellent overview concerning the various approaches invented up to around 2000 was given in [18], published in 2007. 
The historically first algorithm which has a complexity of $O(B^{1+\varepsilon })$ for a search bound of $B$  is the method of $R$. Heath-Brown [6].
 
For $d>2$ Kenji Koyama [7] has generated a large table of integer solutions of  $a^3+b^3+c^3=d$ for noncubes   in the range $1\leqslant d\leqslant 1000$ and $\vert a \vert\leqslant \vert b \vert\leqslant\vert c \vert\leqslant 2^{21}-1$ consists of two tables: Table~1 (55 pages) contains the integer solutions, sorted by $d$, and Table~2 (2 pages) lists the number of primitive solutions found for each $d$ in the search range. 

Consider now some specific: $d=m^3$, $d=m^{12}$ and $d=2m^9$ type values of $d$. Multiply both sides of \eqref{2-3} 
by $m^9$, and apply the change of variable $t\to t/m$ to obtain the more general solution 
\begin{equation}\label{2-5}
(6t^3+m^3)^3-(6t^3-m^3)^3-(6t^2)^3=2m^9
\end{equation}
which is primitive for $\text{GCD}(6t,m)=1$. If $\text{GCD}(6t,m)>1$, then dividing \eqref{2-4} by 
$(\text{GCD}(6t^3, m^3))^3$ gives a primitive solution. For example, for $l, k\geqslant 1$ the solutions
\begin{align}
& (3t^3+2^{3l-1}m^3)^3-(3t^3-2^{3l-1}m^3)3-(2^l 3mt^2)^3=2^{9l-2}m^9 \label{2-6}\\
& (2t^3+3^{3k-1}m^3)^3-(2t^3-3^{3k-1}m^3)^3-(23^k mt^2)^3=23^{9k-3}m^9 \label{2-7}\\
& (t^3+2^{3l-1}3^3k-1 m^3)^3-(t^3-2^{3l-1} 3^{3k-1}m^3)^3-(2^l 3^k km t^2)^3=\notag \\
& \hskip200pt   =2^{9l-2}3^{9k-3}m^9 \label{2-8}
\end{align}
are primitive for $\text{GCD}(3t, 2m)=1$, $\text{GCD}(2t,3m)=1$ and $\text{GCD}(t,6m)=1$ respectively.\\
Equations \eqref{2-5}-\eqref{2-7} give polynomial families for $n = 2, 128, 1458$, $65536, 93312$, $3906250$,  $28697814, \ldots$

An analogous procedure may be applied to \eqref{2-4} to obtain families of solutions for numbers of the form 
$m^{12}$. Multiplying both sides by $m^{12}$ and applying the transformation $t\to t/m$ gives
\begin{equation}\label{2-9}
(9mt^3+m^4)^3-(9t^4+3mt)^3+(9t^4)^3=m^{12}
\end{equation}
which is primitive for $\text{GCD}(3t, m)=1$. In particular, for $3-m$ and $k\geqslant 1$,
\begin{equation}\label{2-10}
(3^k mt^3+3^{4k-2}m^4)^3-(t^4+3^{3k-1}m^3t)^3+(t^4)^3=3^{12k-6}m^{12}
\end{equation}
is primitive for $\text{GCD}(t,3m)=1$. Equations \eqref{2-8} and \eqref{2-9} give families of solutions for
$n = 1, 729, 4096, 2985984, 16777216, 244140625, 387420489, \ldots$

\section{\large New method and results}

In this section a new method and results are surveyed. Here we consider more general framework of the problem of sums of three cubes. We are seeking $P_1(y), P_2(y), P_3(y)$ with the highest possible degree polynomials with integer coefficients  and $Q(y)$ with the lowest possible degree polynomial, so that 
$$P_1^3(y)+P_2^3(y)+P_3^3(y)=Q(y).$$
Actually the solution of this problem has close relation with the above trivial problem, since the case of $\deg Q(y)=0$ coincides with our problem. Nevertheless the estimation of possibility of minimization of $\deg Q(y)$ itself is also an interesting problem.
\vskip3pt
\noindent {\bf RESULT 1.} The first result of this paper is devoted to the case of degrees $(8, 8, 6)$. We search the desired polynomials within the class of polynomials of the form 
$$(a x^8+b x^5+c x^2)^3-(a x^8+b_1 x^5+c_1 x^2)^3-(A x^6+B x^3+C)^3 .$$
\vskip3pt

\noindent
First we expand it
\begin{align*}
&-C^3-3BC^2x^3+(c^3-3B^2C-3AC^2-c1^3) x^6+\\
& +(-B^3+3bc^2 -6ABC-3b1c1^2)x^9 +\\
&+(-3AB^2+3b^2c+3ac^2-3A^2C-3b1^2c1-3ac1^2)x^{12}+\\
&+(b^3-3A^2B-b1^3+6abc-6ab1c1)x^{15}+\\
&+(-A^3+3ab^2-3ab1^2+3a^2c-3a^2c1)x^{18}+(3a^2b-3a^2b1)x^{21}.
\end{align*}
\vskip2pt 
Then we take $b1=b$, $c_1=\frac{-A^3+3a^2c}{3a^2}$, $B=\frac{2Ab}{3a}$ , $C=\frac{-A^4-aAb^2+6a^2Ac}{9a^3}$  and obtain the form
\begin{align*}
& \frac{A^{12}}{729 a^9}+\frac{A^9 b^2}{243 a^8}+\frac{A^6 b^4}{243 a^7}+\frac{A^3 b^6}{729 a^6}-
\frac{2A^9 c}{81 a^7}-\frac{4A^6 b^2 c}{81 a^6} -\frac{2A^3 b^4 c}{81 a^5}+\\
& +\frac{4A^6 c^2}{27a^5} +\frac{4A^3 b^2 c^2}{27a^4}-\frac{8A^3c^3}{27a^3}+\\
& +\left( -\frac{2A^9b}{81a^7}-\frac{4A^6b^3}{81a^6}-\frac{2A^3b^5}{81a^5}+\frac{8A^6bc}{27a^5}+
\frac{8A^3b^3c}{27a^4}-\frac{8A^3bc^2}{9a^3} 
\right) x^3+\\
& +\left(\frac{2A^6b^2}{27a^5}+\frac{A^3b^4}{9a^4}+\frac{A^6 c}{9a^4}-
\frac{4A^3b^2c}{9a^3}-\frac{A^3c^2}{3a^2}\right) x^6 +\\
& +
\left(\frac{A^6b}{9a^4}+\frac{4A^3b^3}{27a^3}-\frac{2A^3bc}{3a^2}\right) x^9
\end{align*}
Further considerations are devoted to the finding of cases interesting for us.
\vskip3pt
\noindent {\it CASE 1:} $b=0$. The result has the form
$$\frac{A^{12}}{729 a^9}-\frac{2A^9 c}{81 a^7}+
\frac{4A^6 c^2}{27a^5}-\frac{8A^3c^3}{27a^3}+\left(
\frac{A^6 c}{9a^4} - \frac{A^3c^2}{3a^2}\right) x^6$$
\vskip3pt
\noindent {\it SUBCASE 1.1:}  $c=0$. The result obtains the form  $\frac{A^{12}}{729a^9}$, which is a cube of an integer, so it is primitive (not interesting).
\vskip3pt
\noindent {\it SUBCASE 1.2:}   $c=\frac{A^3}{3a^2}$,  the result is $-\frac{A^{12}}{729a^9}$, again primitive.
\vskip3pt
\noindent {\it CASE 2:}  $c=\frac{3{{A}^{3}}+4a{{b}^{2}}}{18{{a}^{2}}}$. The result: 
$$-\frac{A^3b^6}{19\,683a^6}-\frac{2A^3b^5x^3}{729a^5}+
\left(\frac{A^9}{108a^6}-\frac{A^3b^4}{243a^4}
\right) x^6$$
\vskip3pt
\noindent Factor the last term, then the result obtains the form
$$- \frac{A^3(-3A^3+2ab^2)(3A^3+2ab^2)}{972a^6}\, .$$
We do the substitution $a=\frac{3A^3}{2b^2}$. Then the result will get the form: 
$$-\frac{64b^{18}}{14\,348\,907A^6}-\frac{64b^{15}x^3}{177\,147A^3}$$
Finally taking $x=\frac{2yb}{A}$ we obtain the result interesting for us:  $-1-648y^3$.
\vskip3pt
\noindent {\it Practical considerations.} Now it’s the time to investigate this result for applications in the solving process of the equation \eqref{2-1}. Since for $\max [abs[a,b,c]]\leqslant 10^{14}$ there are well known tables in [12], so we seek solutions of \eqref{2-1} satisfying the condition $\max [abs[a,b,c]]\geqslant 10^{15}$ with possible small values of 
$abs[d]$ (desirably less than 1000).
 
First rewrite the result:
\begin{align*}
& (54y^2(1+36y^3+432y^6))^3-(18y^2(1+108y^3+1296y^6))^3-\\ 
& -(1+216y^3+3888y^6)^3=-1-648y^3\end{align*}
Of course, calculations expected to be significantly hard, so we'll use Mathematica 11.0 code:
\\
G8$[y\underline{\;\,}]:=$
\\
$1/$GCD$\left[(54y^2(1+36y^3+432y^6)), (18y^2(1+108y^3+1296y^6)), \right.$
\\
$\left.(1+216y^3+3888y^6)\right]$\\
F8$[y\underline{\;\,}]:=$
\\
G8$[y]*\left\{(54y^2(1+36y^3+432y^6)), (18y^2(1+108y^3+1296y^6)), \right.$
\\
$\left.(1+216y^3+3888y^6)\right\}$\\
V8$[y\underline{\;\,}]:=$G8$[y]^{\scriptscriptstyle{\wedge}}3 *(-1-648 y^3)$
\vskip3pt

\noindent For $\{ i=-50, i\leqslant 50, i++,$\\
If$[\text{Abs}[\text{V8}[i]]<1\,000\,000$, If$[\text{Max}[\text{Abs}[\text{F8}[i]]]>1\,000\,000\,000\,000$, \\
Print$[\{i, \text{F8}[i], \text{V8}[i]\}]]]]$
\vskip3pt

\noindent 
The result is:\\
$\{-11, \{5\,000\,250\,899\,358, 5\,000\,250\,895\,002, 6\,887\,541\,673\}, 862\,487\}$\\
$\{-10, \{2\,332\,605\,605\,400, 2\,332\,605\,601\,800, 3\,887\,784\,001\}, 647\,999\}$\\
$\{-9, \{1\,004\,079\,120\,606, 1\,004\,079\,117\,690, 2\,066\,085\,145\}, 472\,391\}$\\
$\{9, \{1\,004\,308\,703\,118, 1\,004\,308\,700\,202, 2\,066\,400\,073\}, -472\,393\}$\\
$\{10, \{2\,332\,994\,405\,400, 2\,332\,994\,401\,800, 3\,888\,216\,001\}, -648\,001\}$\\
$\{11, \{5\,000\,877\,065\,646, 5\,000\,877\,061\,290, 6\,888\,116\,665\}, -862\,489\}$

\noindent This means that,  for example
$$ 1004079120606^{\scriptscriptstyle{\wedge}}3 - 1004079117690^{\scriptscriptstyle{\wedge}}3- 2066085145^{\scriptscriptstyle{\wedge}}3 = 472391$$

\noindent {\bf RESULT 2.}
More enhanced result is obtained for the case $(9, 9, 7)$:\\
Expand$[(3(1+120 y^3 +3456y^6+31\,104 y^9))^{\scriptscriptstyle{\wedge}}3-
((-1+216y^3+10\,368 y^6+93\,312y^9))^{\scriptscriptstyle{\wedge}}3 -$

$(4y(5+324y^3+3888y^6))^{\scriptscriptstyle{\wedge}}3]$\\
$28+1072y^3$\\
G9$[y\underline{\;\,}]:=$
\\
$1/$GCD$[(3(1+120y^3+3456y^6+31\,104y^9)), ((-1+216y^3+10\,368 y^6+93\,312y^9)),$
\\
$(4y(5+324y^3+3888y^6))]$\\
F9$[y\underline{\;\,}]:=$
\\
G9$[y]*\{(3(1+120y^3+3456y^6+31\,104y^9)), ((-1+216y^3+10\,368 y^6+93\,312y^9)),$
\\
$(4y(5+324y^3+3888y^6))\}$\\
V9$[y\underline{\;\,}]:=$
G9$[y]^{\scriptscriptstyle{\wedge}}3 *(28+1072y^3)$
\vskip3pt

\noindent For $\{ i=-50, i\leqslant 50, i++,$\\
If$[\text{Abs}[\text{V9}[i]]<10\,000\,000$, If$[\text{Max}[\text{Abs}[\text{F9}[i]]]>10\,000\,000\,000\,000$, \\
Print$[\{i, \text{F9}[i], \text{V9}[i]\}]]]]$
\vskip3pt

\noindent 
The result is:\\
$\{-21, \{-74\,114\,970\,486\,757\,701, -74\,114\,970\,485\,424\,121, -28\,010\,276\,942\,676\},$

$ -9\,927\,764\}$\\
$\{-20, \{-47\,775\,080\,450\,879\,997, -47\,775\,080\,449\,728\,001, -19\,906\,352\,640\,400\}, $

$-8\,575\,972\}$\\
$\{-19, \{-30\,110\,146\,685\,929\,077, -30\,110\,146\,684\,941\,385, -13\,901\,324\,389\,292\}, $

$-7\,352\,820\}$\\
$\{-18, \{-18\,508\,949\,466\,179\,901, -18\,508\,949\,465\,340\,097, -9\,521\,109\,889\,128\}, $

$-6\,251\,876\}$\\
$\{-17, \{-11\,065\,421\,675\,141\,349, -11\,065\,421\,674\,433\,881, -6\,381\,478\,799\,620\}, $

$-5\,266\,708\}$\\
$\{-16, \{-6\,412\,177\,868\,488\,701, -6\,412\,177\,867\,898\,881, -4\,174\,623\,277\,376\}, $

$-4\,390\,884\}$\\
$\{-15, \{-3\,587\,108\,653\,214\,997, -3\,587\,108\,652\,729\,001, -2\,657\,139\,390\,300\},$

$ -3\,617\,972\}$\\
$\{-14, \{-1\,927\,845\,532\,267\,197, -1\,927\,845\,531\,872\,065, -1\,639\,341\,027\,352\}, $

$-2\,941\,540\}$\\
$\{14, \{1\,928\,001\,664\,725\,699, 1\, 928\,001\,664\,330\,559, 1\,639\,440\,601\,624\}, 2\,941\,596\}$\\
$\{15, \{3\,587\,344\,849\,215\,003, 3\,587\,344\,848\,728\,999, 2\,657\,270\,610\,300\}, 3\,618\,028\}$\\
$\{16, \{6\,412\,525\,760\,839\,683, 6\,412\,525\,760\,249\,855, 4\,174\,793\,146\,688\}, 4\,390\,940\}$\\
$\{17, \{11\,065\,922\,191\,772\,139, 11\,065\,922\,191\,064\,663, 6\,381\,695\,286\,052\}, 5\,266\,764\}$\\
$\{18, \{18\,509\,654\,743\,656\,771, 18\,509\,654\,742\,816\,959, 9\,521\,381\,986\,920\}, 6\,251\,932\}$\\
$\{19, \{30\,111\,122\,229\,317\,499, 30\,111\,122\,228\,329\,799, 13\,901\,662\,181, 324\}, 7\,352\,876\}$\\
$\{20, \{47\,776\,407\,554\,880\,003, 47\,776\,407\,553\,727\,999, 19\,906\,767\,360\,400\}, 8\,576\,028\}$\\
$\{21, \{74\,116\,748\,933\,042\,763, 74\,116\,748\,931\,709\,175, 28\,010\,781\,037\,428\}, 9\,927\,820\}$

\noindent Here the most interesting triple is:
$$1928001664725699^{\scriptscriptstyle{\wedge}}3 -1928001664330559^{\scriptscriptstyle{\wedge}}3 - 1639440601624^{\scriptscriptstyle{\wedge}}3 = 2941596$$
\noindent {\bf RESULT 3.} Consider now the case $(25, 25, 18)$. Through the same way we get:\\
Expand$\left[\left(\frac{1}{18}y(63+36 y^3 +280y^6+672y^{12}+768y^{18}+512y^{24})\right)^{\scriptscriptstyle{\wedge}}3-
\right.$
\\
$\left(\frac{1}{18}y(63-36 y^3 +280y^6+672y^{12}+768y^{18}+512y^{24})\right)^{\scriptscriptstyle{\wedge}}3-$
\\
$\left.\left(\frac{2}{3}(3+20y^6+32y^{12}+32y^{18})\right)^{\scriptscriptstyle{\wedge}}3\right]$
\\
$-8-13y^6$\\
Then we use the code:\\
G25$[y\underline{\;\,}]:=$
\\
$1/$GCD$\left[\left(\frac{1}{18}y(63+36y^3+280y^6+672y^{12}+768y^{18}+512y^{24})\right),\right.$
\\
$\left(\frac{1}{18}y(63-36 y^3 +280y^6+672y^{12}+768y^{18}+512y^{24})\right),$
\\
$\left.\left(\frac{2}{3}(3+20y^6+32y^{12}+32y^{18})\right)\right]$
\\
F25$[y\underline{\;\,}]:=$
\\
G25$[y]*\left\{\left(\frac{1}{18}y(63+36 y^3 +280y^6+672y^{12}+768y^{18}+512y^{24})\right) \right.$
\\
$\left(
\frac{1}{18}y(63-36 y^3 +280y^6+672y^{12}+768y^{18}+512y^{24})\right),$
\\
$\left.\left(\frac{2}{3}(3+20y^6+32y^{12}+32y^{18})\right)\right\}$\\
V25$[y\underline{\;\,}]:=$
G25$[y]^{\scriptscriptstyle{\wedge}}3 *(-8-13y^6)$
\vskip3pt

\noindent For $\{ i=-50, i\leqslant 0, i++$, If$[\text{Abs}[\text{V25}[i]]<1\,000\,000\,000$,\\ 
If$[\text{Max}[\text{Abs}[\text{F25}[i]]]>100\,000\,000\,000\,000\,000\,000$, \\
Print$[\{i, \text{F25}[i], \text{V25}[i]\}]]]]$
\vskip3pt

\noindent 
The result is:\\
$\{-28,\{-21\,474\,261\,883\,010\,575\,951\,072\,188\,890\,073\,079\,601,$

$-21\,474\,261\,883\,010\,575\,951\,072\,188\,890\,074\,308\,913,$

$1\,193\,639\,792\,964\,388\,519\,010\,222\,081\}, -783\,071\,745\}$\\
$\{-24, \{-455\,248\,482\,071\,553\,635\,586\,938\,761\,943\,642\,154,$

$-455\,248\,482\,071\,553\,635\,586\,938\,761\,944\,305\,706,$

$74\,444\,235\,905\,117\,108\,623\,638\,529\}, -310\,542\,337\}$\\
$\{-20, \{-4\,772\,185\,996\,292\,552\,640\,284\,454\,399\,840\,035,$

$-4\,772\,185\,996\,292\,552\,640\,284\,454\,400\,160\,035,$

$2\,796\,202\,710\,357\,333\,760\,000\,001\}, -104\,000\,001\}$\\
$\{-18, \{-685\,185\,558\,884\,868\,266\,867\,584\,546\,812\,447,$

$-685\,188\,558\,884\,868\,266\,867\,584\,547\,232\,351,$

$839\,390\,063\,618\,729\,392\,197\,122\}, -442\,158\,920\}$\\
$\{-16, \{-18\,028\,810\,148\,480\,439\,485\,764\,921\,655\,324,$

$-18\,028\,810\,148\,480\,439\,485\,764\,921\,786\,396, 50\,371\,912\,153\,009\,412\,374\,529\},$

$ -27\,262\,977\}$\\
$\{-14, \{-1\,279\,966\,002\,355\,352\,271\,319\,733\,014\,033,$

$-1\,279\,966\,002\,355\,352\,271\,319\,733\,167\,697, 9\,106\,750\,099\,297\,148\,243\,458\},$

$ -97\,883\,976\}$\\
$\{-13, \{-401\,431\,450\,040\,755\,185\,937\,727\,808\,239,$

$-401\,431\,450\,040\,755\,185\,937\,728\,036\,727, 4\,798\,098\,357\,335\,276\,047\,604\}, $

$-501\,988\,200\}$\\
$\{-12, \{-13\,567\,468\,738\,287\,354\,512\,175\,542\,037,$

$-13\,567\,468\,738\,287\,354\,512\,175\,583\,509, 283\,982\,316\,767\,867\,289\,601\},$

$ -4\,852\,225\}$\\
$\{-11, \{-6\,163\,749\,044\,483\,681\,037\,311\,693\,681,$

$-6\,163\,749\,044\,483\,681\,037\,311\,810\,809, 237\,223\,272\,615\,326\,524\,916\},$

$ -184\,242\,408\}$\\
$\{-10, \{-284\,444\,871\,111\,484\,444\,599\,980\,035,$

$-284\,444\,871\,111\,484\,444\,600\,020\,035, 21\,333\,354\,666\,680\,000\,002\},$

$ -13\,000\,008\}$\\
$\{-9, \{-40\,840\,534\,128\,228\,425\,942\,658\,291,$

$-40\,840\,534\,128\,228\,425\,942\,710\,779, 6\,404\,049\,823\,013\,024\,116\},$

$ -55\,269\,928\}$\\
$\{-8, \{-537\,303\,438\,740\,776\,265\,183\,246,$

$-537\,303\,438\,740\,776\,265\,191\,438, 192\,154\,317\,110\,640\,641\},$
$ -425\,985\}$\\
$\{-7, \{-76\,292\,876\,409\,395\,782\,365\,173,$

$-76\,292\,876\,409\,395\,782\,384\,381,$
$69\,479\,570\,742\,237\,044\}, -12\,235\,560\}$\\
$\{-6, \{-808\,709\,748\,243\,399\,993\,333,$

$-808\,709\,748\,243\,399\,998\,517,$
$2\,166\,658\,847\,571\,458\}, -606\,536\}$

\vskip3pt \noindent {\bf RESULT 4.} Case $(27, 27, 20)$. \\
Expand$\left[\right.$
\\
$(1+177\,710\,598y^3+17\,738\,799\,316\,992y^6+7\,466\,750\,649\,114\,265\,387\,008y^9+$
\\
$123\,267\,709\,616\,967\,231\,382\,892\,912\,798\,859\,264y^{15}+$
\\ 
$945\,048\,866\,667\,847\,329\,755\,658\,073\,857\,921\,409\,357\,870\,792\,704 y^{21}+$
\\
$2\,867\,531\,822\,071\,470\,533\,473\,102\,364\,531\,302\,968\,788\,957\,393\,604\,240\,929\,193\,984 
y^{27})^{\scriptscriptstyle{\wedge}}3-$
\\
$(-1+177\,710\,598 y^3-17\,738\,799\,316\,992 y^6 +
7\,466\,750\,649\,114\,265\,387\,008 y^9+$
\\
$123\,267\,709\,616\,967\,231\,382\,892\,912\,798\,859\,264y^{15}+$
\\
$945\,048\,866\,667\,847\,329\,755\,658\,073\,857\,921\,409\,357\,870\,792\,704 y^{21}+$
\\
$2\,867\,531\,822\,071\,470\,533\,473\,102\,364\,531\,302\,968\,788\,957\,393\,604\,240\,929\,193\,984 
y^{27})^{\scriptscriptstyle{\wedge}}3-$
\\
$(234y^2(2455+83\,248\,185\,194\,643\,456 y^6+
974\,937\,062\,077\,718\,261\,926\,943\,784\,960 y^{12}+$
\\
$4\,087\,723\,046\,938\,680\,100\,330\,712\,454\,833\,737\,367\,027\,712 y^{18}))^{\scriptscriptstyle{\wedge}}3]$
\\
$2+8\,606\,991\,384\,576 y^6$
\vskip3pt

\noindent G27$[y\underline{\;\,}]:=$
\\
$1/$GCD$[(1+177\,710\,598y^3+17\,738\,799\,316\,992y^6+7\,466\,750\,649\,114\,265\,387\,008y^9+$
\\
$123\,267\,709\,616\,967\,231\,382\,892\,912\,798\,859\,264y^{15}+$
\\ 
$945\,048\,866\,667\,847\,329\,755\,658\,073\,857\,921\,409\,357\,870\,792\,704 y^{21}+$
\\
$2\,867\,531\,822\,071\,470\,533\,473\,102\,364\,531\,302\,968\,788\,957\,393\,604\,240\,929\,193\,984 
y^{27})$,
\\
$(-1+177\,710\,598 y^3-17\,738\,799\,316\,992 y^6 +
7\,466\,750\,649\,114\,265\,387\,008 y^9+$
\\
$123\,267\,709\,616\,967\,231\,382\,892\,912\,798\,859\,264y^{15}+$
\\ 
$945\,048\,866\,667\,847\,329\,755\,658\,073\,857\,921\,409\,357\,870\,792\,704 y^{21}+$
\\
$2\,867\,531\,822\,071\,470\,533\,473\,102\,364\,531\,302\,968\,788\,957\,393\,604\,240\,929\,193\,984 
y^{27})$,
\\
$(234y^2(2455+83\,248\,185\,194\,643\,456 y^6+
974\,937\,062\,077\,718\,261\,926\,943\,784\,960 y^{12}+$
\\
$4\,087\,723\,046\,938\,680\,100\,330\,712\,454\,833\,737\,367\,027\,712 y^{18}))]$
\\
F27$[y\underline{\;\,}]:=$
\\
G27$[y]*\{(1+177\,710\,598y^3+17\,738\,799\,316\,992y^6+7\,466\,750\,649\,114\,265\,387\,008y^9+$
\\
$123\,267\,709\,616\,967\,231\,382\,892\,912\,798\,859\,264y^{15}+$
\\
$945\,048\,866\,667\,847\,329\,755\,658\,073\,857\,921\,409\,357\,870\,792\,704 y^{21}+$
\\
$2\,867\,531\,822\,071\,470\,533\,473\,102\,364\,531\,302\,968\,788\,957\,393\,604\,240\,929\,193\,984 
y^{27})$,
\\
$(-1+177\,710\,598 y^3-17\,738\,799\,316\,992 y^6 +
7\,466\,750\,649\,114\,265\,387\,008 y^9+$
\\
$123\,267\,709\,616\,967\,231\,382\,892\,912\,798\,859\,264y^{15}+$
\\
$945\,048\,866\,667\,847\,329\,755\,658\,073\,857\,921\,409\,357\,870\,792\,704 y^{21}+$
\\
$2\,867\,531\,822\,071\,470\,533\,473\,102\,364\,531\,302\,968\,788\,957\,393\,604\,240\,929\,193\,984 
y^{27})$,
\\
$(234y^2(2455+83\,248\,185\,194\,643\,456 y^6+
974\,937\,062\,077\,718\,261\,926\,943\,784\,960 y^{12}+$
\\
$4\,087\,723\,046\,938\,680\,100\,330\,712\,454\,833\,737\,367\,027\,712 y^{18}))\}$
\\
V27$[y\underline{\;\,}]:=$
G27$[y]^{\scriptscriptstyle{\wedge}}3 *(2+8\,606\,991\,384\,576 y^6)$
\vskip3pt

\noindent For $\{ i=1, i\leqslant 10, i++$,
\\
For $\{ j=1, j\leqslant 100, j++$,
\\
If$[\text{GCD}[i, j]==1,$ If$[\text{Abs}[\text{V27}[i/j]]<10\,000\,000\,000\,000$,
\\
If$[\text{Max}[\text{Abs}[\text{F27}[i/j]]]>100\,000\,000\,000\,000\,000\,000$, 
\\
Print$[\{i/j, \text{F27}[i/j], \text{V27}[i/j]\}]]]]]]$
\vskip3pt

\noindent 
The result:\\ 
$\left\{1, \{2\,867\,531\,822\,072\,415\,582\,339\,770\,335\,128\,768\,243\,837\,573\,448\,573\,364\,841\,260\,551,\right.$

$2\,867\,531\,822\,072\,415\,582\,339\,770\,335\,128\,768\,243\,837\,513\,448\,537\,887\,242\,626\,565, $

$\left. 956\,527\,192\,983\,879\,278\,749\,912\,919\,984\,460\,784\,277\,308\,422\}, 8\,606\,991\,384\,578\right\}$\\
$\left\{\frac{1}{2}, \{85\,459\,107\,820\,151\,855\,377\,565\,309\,350\,915\,495\,772\,943\,736\,866\,829\,063,\right.$

$85\,459\,107\,820\,151\,855\,377\,565\,309\,350\,915\,495\,772\,941\,519\,516\,914\,431, $

$\left. 3\,648\,861\,667\,626\,387\,790\,371\,484\,426\,537\,654\,889\,478\}, 8\,606\,991\,384\,704\right\}$\\
$\left\{\frac{1}{3}, \{3\,384\,362\,531\,021\,750\,766\,977\,249\,652\,864\,977\,396\,620\,327\,182\,859, \right.$

$3\,384\,362\,531\,021\,750\,766\,977\,249\,652\,864\,977\,396\,182\,332\,137\,977, $

$\left. 2\,468\,963\,878\,546\,861\,668\,845\,110\,751\,426\,429\,958\}, 8\,606\,991\,386\,034\right\}$\\
$\left\{\frac{1}{6}, \{ 100\,861\,864\,753\,002\,334\,438\,937\,135\,077\,135\,856\,775\,461 ,\right.$

$ 100\,861\,864\,753\,002\,334\,438\,937\,135\,077\,108\,482\,085\,085 ,$

$\left. 9\,418\,349\,859\,585\,944\,691\,803\,620\,862\,982\}, 8\,606\,991\,477\,888\right\}$\\
$\left\{\frac{1}{13}, \{31\,259\,481\,996\,382\,783\,282\,485\,723\,728\,803, \right.$

$31\,259\,481\,996\,382\,783\,282\,485\,628\,177\,289, $
$654\,295\,791\,943\,512\,069\,680\,286\}, $

$\left. 3\,917\,615\,402\right\}$\\
$\left\{\frac{1}{26}, \{931\,698\,527\,611\,284\,677\,846\,527,\right.$

$931\,698\,527\,611\,284\,671\,874\,455,$
$\left. 2\,496\,121\,394\,505\,512\,094 \}, 3\,917\,892\,224\right\}$\\
$\left\{\frac{2}{13}, \{4\,195\,570\,082\,123\,558\,929\,847\,343\,917\,618\,224\,454\,845,\right.$

$4\,195\,570\,082\,123\,558\,929\,847\,343\,917\,612\,109\,159\,587, $

$\left. 686\,078\,086\,852\,418\,912\,046\,595\,355\,256\}, 250\,727\,108\,906\right\}$\\
$\left\{\frac{2}{39}, \{4\,951\,849\,882\,748\,149\,562\,576\,914\,949,\right.$
$ 4\,951\,849\,882\,748\,149\,562\,501\,417\,243,$

$\left. 1\,770\,910\,652\,870\,714\,110\,584\}, 250\,730\,307\,738\right\}$\\
$\left\{\frac{3}{13}, \{238\,371\,848\,619\,844\,446\,597\,099\,170\,950\,280\,787\,231\,652\,063,\right.$

$ 238\,371\,848\,619\,844\,446\,597\,099\,170\,950\,280\,717\,574\,617\,285,$

$\left. 2\,281\,385\,738\,232\,016\,715\,434\,029\,491\,949\,966\}, 2\,855\,938\,429\,226\right\}$\\
$\left\{\frac{3}{26}, \{7\,104\,035\,657\,336\,443\,018\,251\,070\,280\,549\,440\,669,\right.$

$7\,104\,035\,657\,336\,443\,018\,251\,070\,276\,195\,875\,893, $

$\left. 8\,702\,796\,803\,288\,457\,118\,679\,259\,534\}, 2\,855\,938\,706\,048\right\}$

\section{\large Summary}
 
Denote $m={\max \{\deg (P_1), \deg (P_2), \deg (P_3)\}}/{\deg (Q})$,  where $P_{1}^{3}(y)+P_{2}^{3}(y)+P_{3}^{3}(y)=Q(y)$ and $n=\text{IntegerPart}\left[\frac{\text{Max}[\text{Abs}[a, b, c]]}{d}\right]$,  where 
$a^3+b^3+c^3=d$ as a describers of solution quality. Note that if $\sup n(d)=+\infty $, then the equation 
$a^3+b^3+c^3=d$ has infinite set of solutions.

We consider the values of  $m$ and $n$ for each result. It is clear, that as much greater are $n$ and  $m$ as cool is result.
\vskip3pt
\noindent {\it Result 1.} We have
$$(54y^2+1944y^5+23328y^8)^3-(18y^2+1944y^5+23328y^8)^3-(1+216y^3+3888y^6)^3=-1-648y^3$$  
where $m=\frac{8}{3}$. Consider the code:\\
G8$[y\underline{\;\,}]:=$
\\
$1/$GCD$\left[(54y^2(1+36y^3+432y^6)), (18y^2(1+108y^3+1296y^6)), (1+216y^3+3888y^6)\right]$\\
F8$[y\underline{\;\,}]:=$
\\
G8$[y]*\left\{(54y^2(1+36y^3+432y^6)), (18y^2(1+108y^3+1296y^6)), (1+216y^3+3888y^6)\right\}$\\
V8$[y\underline{\;\,}]:=$G8$[y]^{\scriptscriptstyle{\wedge}}3 *(-1-648 y^3)$

\vskip3pt \noindent For $\{ i=-100, i\leqslant 100, i++,$\\
If$[\text{Abs}[\text{V8}[i]]<10\,000$, If$[\text{Max}[\text{Abs}[\text{F8}[i]]]>1\,000\,000$, 
Print$[\{i, \text{F8}[i], \text{V8}[i]\}]]]]$
\vskip3pt \noindent 
The result is:\\
$\{-2, \{5\,909\,976, 5\,909\,832, 247\,105\}, 5183\}$\\
$\{2, \{6\,034\,392, 6\,034\,248, 250\,561\}, -5185\}$

\noindent So $n=\text{IntegerPart}\left[ \frac{5909976}{5183} \right]=1140$.
\vskip3pt
\noindent {\it Result 2.} We have
\begin{align*}
&(3+360y^3+10368y^6+93312y^9)^3-(-1+216y^3+10368y^6+93312y^9)^3- \\ 
&-(20y+1296y^4+15552y^7)^3=28+1072y^3  
\end{align*}
where $m=3$. Now consider the code:\\
G9$[y\underline{\;\,}]:=$
\\
$1/$GCD$[(3(1+120y^3+3456y^6+31\,104y^9)), ((-1+216y^3+10\,368 y^6+93\,312y^9)),$
\\
$(4y(5+324y^3+3888y^6))]$\\
F9$[y\underline{\;\,}]:=$
\\
G9$[y]*\{(3(1+120y^3+3456y^6+31\,104y^9)), ((-1+216y^3+10\,368 y^6+93\,312y^9)),$
\\
$(4y(5+324y^3+3888y^6))\}$\\
V9$[y\underline{\;\,}]:=$
G9$[y]^{\scriptscriptstyle{\wedge}}3 *(28+1072y^3)$
\vskip3pt

\noindent For $\{ i=-50, i\leqslant 50, i++,$\\
If$[\text{Abs}[\text{V9}[i]]<10\,000$, If$[\text{Max}[\text{Abs}[\text{F9}[i]]]>10\,000\,000$, 
Print$[\{i, \text{F9}[i], \text{V9}[i]\}]]]]$

\vskip3pt \noindent 
The result is:\\
$\{-2, \{-47\,115\,069, -47\,113\,921, -1\, 969\,960\}, -8548\}$\\
$\{2, \{48\,442\,179, 48\,441\,023, 2\,011\,432\}, 8604\}$

\vskip3pt \noindent 
So $n=\text{IntegerPart}\left[ \frac{47115069}{8548} \right]=5511$.
\vskip3pt
noindent {\it Result 3.} We have 
\begin{align*}
&\left(\frac{1}{18}y(63+36y^3+280y^6+672y^{12}+768y^{18}+512y^{24})\right)^3-  
\left(\frac{1}{18}y(63-36y^3+280y^6+\right.\\
&+672y^{12}+768y^{18}+512y^{24})\Big)^3-
\left(\frac{2}{3}(3+20y^6+32y^{12}+32y^{18})\right)^3= -8-13y^6 
\end{align*}
where $m=\frac{25}{6}$.  Note that here the leading coefficient of $Q$ is significantly small. Consider the code:\\
G25$[y\underline{\;\,}]:=$
\\
$1/$GCD$\left[\left(\frac{1}{18}y(63+36y^3+280y^6+672y^{12}+768y^{18}+512y^{24})\right),\right.$
\\
$\left.\left(
\frac{1}{18}y(63-36 y^3 +280y^6+672y^{12}+768y^{18}+512y^{24})\right),
\left(\frac{2}{3}(3+20y^6+32y^{12}+32y^{18})\right)\right]$
\\
F25$[y\underline{\;\,}]:=$
\\
G25$[y]*\left\{\left(\frac{1}{18}y(63+36 y^3 +280y^6+672y^{12}+768y^{18}+512y^{24})\right) \right.$
\\
$\left.\left(
\frac{1}{18}y(63-36 y^3 +280y^6+672y^{12}+768y^{18}+512y^{24})\right),
\left(\frac{2}{3}(3+20y^6+32y^{12}+32y^{18})\right)\right\}$\\
V25$[y\underline{\;\,}]:=$
G25$[y]^{\scriptscriptstyle{\wedge}}3 *(-8-13y^6)$
\vskip3pt

\noindent For $\{ i=1, i\leqslant 50, i++$, \\
If$[\text{Abs}[\text{V25}[i]]<1\,000\,000$, 
If$[\text{Max}[\text{Abs}[\text{F25}[i]]]>1\,000\,000\,000\,000\,000$, \\
Print$[\{i, \text{F25}[i], \text{V25}[i]\}]]]]$
\vskip3pt
\noindent
The result is:\\
$\{4. \{16\,018\,663\,989\,936\,391, 16\,018\,663\,989\,935\,879, 733\,186\,736\,129\}, -6657\}$\\
$\{6, \{808\,709\,748\,243\,399\,998\,517, 808\,709\,748\,243\,399\,993\,333, 2\,166\,658\,847\,571\,458\}, -606\,536\}$\\
$\{8, \{537\,303\,438\,740\,776\,265\,191\,438, 537\,303\,438\,740\,776\,265\,183\,246,$

$ 192\,154\,317\,110\,640\,641\}, -425\,985\}$
\vskip3pt
\noindent
So $n=\text{IntegerPart}\left[ \frac{537303438740776265191438}{425985} \right]=1261320090474491508$.
\vskip3pt
\noindent {\it Result 4.} We have
\begin{align*}
&(1+177710598y^3+1773879931\,6992y^6+7466750649\,1142653870\,08y^9+ \\ 
&1232677096\,1696723138\,2892912798\,859264y^{15}+ \\ 
&9450488666\,6784732975\,5658073857\,9214093578\,70792704y^{21}+ \\ 
&2867531822\,0714705334\,7310236453\,1302968788\,9573936042\,4092919398\,4y^{27})^3- \\ 
&-(1+177710598y^3-1773879931\,6992y^6+7466750649\,1142653870\,08y^9+ \\
&1232677096\,1696723138\,2892912798\,859264y^{15}+ \\ 
&9450488666\,6784732975\,5658073857\,9214093578\,70792704y^{21}+ \\ 
&2867531822\,0714705334\,7310236453\,1302968788\,9573936042\,4092919398\,4y^{27})^3- \\
&-(574470y^2+1948007533\,5546568704y^8+2281352725\,2618607329\,0904845680\,640y^{14}+  \\
&+9565271929\,8365114347\,7386714431\,0945438844\,84608y^{20})^3=2+8606991384\,576y^6  
\end{align*}
where $m=\frac{9}{2}$ . Consider the code:
\vskip3pt

\noindent G27$[y\underline{\;\,}]:=$
\\
$1/$GCD$[(1+177\,710\,598y^3+17\,738\,799\,316\,992y^6+7\,466\,750\,649\,114\,265\,387\,008y^9+$
\\
$123\,267\,709\,616\,967\,231\,382\,892\,912\,798\,859\,264y^{15}+$
\\
$945\,048\,866\,667\,847\,329\,755\,658\,073\,857\,921\,409\,357\,870\,792\,704 y^{21}+$
\\
$2\,867\,531\,822\,071\,470\,533\,473\,102\,364\,531\,302\,968\,788\,957\,393\,604\,240\,929\,193\,984 
y^{27})$,
\\
$(-1+177\,710\,598 y^3-17\,738\,799\,316\,992 y^6 +
7\,466\,750\,649\,114\,265\,387\,008 y^9+$
\\
$123\,267\,709\,616\,967\,231\,382\,892\,912\,798\,859\,264y^{15}+$
\\
$945\,048\,866\,667\,847\,329\,755\,658\,073\,857\,921\,409\,357\,870\,792\,704 y^{21}+$
\\
$2\,867\,531\,822\,071\,470\,533\,473\,102\,364\,531\,302\,968\,788\,957\,393\,604\,240\,929\,193\,984 
y^{27})$,

$(234y^2(2455+83\,248\,185\,194\,643\,456 y^6+
974\,937\,062\,077\,718\,261\,926\,943\,784\,960 y^{12}+$
\\
$4\,087\,723\,046\,938\,680\,100\,330\,712\,454\,833\,737\,367\,027\,712 y^{18}))]$
\\
F27$[y\underline{\;\,}]:=$
\\
G27$[y]*\{(1+177\,710\,598y^3+17\,738\,799\,316\,992y^6+7\,466\,750\,649\,114\,265\,387\,008y^9+$
\\
$123\,267\,709\,616\,967\,231\,382\,892\,912\,798\,859\,264y^{15}+$
\\
$945\,048\,866\,667\,847\,329\,755\,658\,073\,857\,921\,409\,357\,870\,792\,704 y^{21}+$
\\
$2\,867\,531\,822\,071\,470\,533\,473\,102\,364\,531\,302\,968\,788\,957\,393\,604\,240\,929\,193\,984 
y^{27})$,
\\
$(-1+177\,710\,598 y^3-17\,738\,799\,316\,992 y^6 +
7\,466\,750\,649\,114\,265\,387\,008 y^9+$
\\
$123\,267\,709\,616\,967\,231\,382\,892\,912\,798\,859\,264y^{15}+$
\\
$945\,048\,866\,667\,847\,329\,755\,658\,073\,857\,921\,409\,357\,870\,792\,704 y^{21}+$
\\
$2\,867\,531\,822\,071\,470\,533\,473\,102\,364\,531\,302\,968\,788\,957\,393\,604\,240\,929\,193\,984 
y^{27})$,
\\
$(234y^2(2455+83\,248\,185\,194\,643\,456 y^6+
974\,937\,062\,077\,718\,261\,926\,943\,784\,960 y^{12}+$
\\
$4\,087\,723\,046\,938\,680\,100\,330\,712\,454\,833\,737\,367\,027\,712 y^{18}))\}$
\\
V27$[y\underline{\;\,}]:=$
G27$[y]^{\scriptscriptstyle{\wedge}}3 *(2+8\,606\,991\,384\,576 y^6)$
\vskip3pt

\noindent For $\{ i=1, i\leqslant 10, i++$,

For $\{ j=1, j\leqslant 100, j++$,
\\
If$[\text{GCD}[i, j]==1,$ If$[\text{Abs}[\text{V27}[i/j]]<10\,000\,000\,000\,000$,
\\
If$[\text{Max}[\text{Abs}[\text{F27}[i/j]]]>100\,000\,000\,000\,000\,000\,000$, 
\\
Print$[\{i/j, \text{F27}[i/j], \text{V27}[i/j]\}]]]]]]$
\vskip3pt

\noindent 
The result:\\  
$\left\{1, \{2\,867\,531\,822\,072\,415\,582\,339\,770\,335\,128\,768\,243\,837\,573\,448\,573\,364\,841\,260\,551,\right.$

$2\,867\,531\,822\,072\,415\,582\,339\,770\,335\,128\,768\,243\,837\,513\,448\,537\,887\,242\,626\,565, $

$\left. 956\,527\,192\,983\,879\,278\,749\,912\,919\,984\,460\,784\,277\,308\,422\}, 8\,606\,991\,384\,578\right\}$\\
$\left\{\frac{1}{2}, \{85\,459\,107\,820\,151\,855\,377\,565\,309\,350\,915\,495\,772\,943\,736\,866\,829\,063,\right.$

$85\,459\,107\,820\,151\,855\,377\,565\,309\,350\,915\,495\,772\,941\,519\,516\,914\,431, $

$\left. 3\,648\,861\,667\,626\,387\,790\,371\,484\,426\,537\,654\,889\,478\}, 8\,606\,991\,384\,704\right\}$\\
$\left\{\frac{1}{3}, \{3\,384\,362\,531\,021\,750\,766\,977\,249\,652\,864\,977\,396\,620\,327\,182\,859, \right.$

$3\,384\,362\,531\,021\,750\,766\,977\,249\,652\,864\,977\,396\,182\,332\,137\,977, $

$\left. 2\,468\,963\,878\,546\,861\,668\,845\,110\,751\,426\,429\,958\}, 8\,606\,991\,386\,034\right\}$\\
$\left\{\frac{1}{6}, \{ 100\,861\,864\,753\,002\,334\,438\,937\,135\,077\,135\,856\,775\,461 ,\right.$

$ 100\,861\,864\,753\,002\,334\,438\,937\,135\,077\,108\,482\,085\,085 ,$

$\left. 9\,418\,349\,859\,585\,944\,691\,803\,620\,862\,982\}, 8\,606\,991\,477\,888\right\}$\\
$\left\{\frac{1}{13}, \{31\,259\,481\,996\,382\,783\,282\,485\,723\,728\,803, \right.$

$31\,259\,481\,996\,382\,783\,282\,485\,628\,177\,289, $
$\left. 654\,295\,791\,943\,512\,069\,680\,286\}, 3\,917\,615\,402\right\}$\\
$\left\{\frac{1}{26}, \{931\,698\,527\,611\,284\,677\,846\,527,\right.$

$931\,698\,527\,611\,284\,671\,874\,455,$
$\left. 2\,496\,121\,394\,505\,512\,094 \}, 3\,917\,892\,224\right\}$\\
$\left\{\frac{2}{13}, \{4\,195\,570\,082\,123\,558\,929\,847\,343\,917\,618\,224\,454\,845,\right.$

$4\,195\,570\,082\,123\,558\,929\,847\,343\,917\,612\,109\,159\,587, $

$\left. 686\,078\,086\,852\,418\,912\,046\,595\,355\,256\}, 250\,727\,108\,906\right\}$\\
$\left\{\frac{2}{39}, \{4\,951\,849\,882\,748\,149\,562\,576\,914\,949,\right.$

$ 4\,951\,849\,882\,748\,149\,562\,501\,417\,243,$
$\left. 1\,770\,910\,652\,870\,714\,110\,584\}, 250\,730\,307\,738\right\}$\\
$\left\{\frac{3}{13}, \{238\,371\,848\,619\,844\,446\,597\,099\,170\,950\,280\,787\,231\,652\,063,\right.$

$ 238\,371\,848\,619\,844\,446\,597\,099\,170\,950\,280\,717\,574\,617\,285,$

$\left. 2\,281\,385\,738\,232\,016\,715\,434\,029\,491\,949\,966\}, 2\,855\,938\,429\,226\right\}$\\
$\left\{\frac{3}{26}, \{7\,104\,035\,657\,336\,443\,018\,251\,070\,280\,549\,440\,669,\right.$

$7\,104\,035\,657\,336\,443\,018\,251\,070\,276\,195\,875\,893, $

$\left. 8\,702\,796\,803\,288\,457\,118\,679\,259\,534\}, 2\,855\,938\,706\,048\right\}$
 
\noindent Note, that nevertheless $8606991384578$ is rather greater than the considered range, however three cubes are essentially big numbers:\\
$2\,867\,531\,822\,072\,415\,582\,339\,770\,335\,128\,768\,243\,837\,513\,448\,573\,364\,841\,260\,551^{\scriptscriptstyle{\wedge}}3 -$\\
$2\,867\,531\,822\,072\,415\,582\,339\,770\,335\,128\,768\,243\,837\,513\,448\,537\,887\,242\,626\,565^{\scriptscriptstyle{\wedge}}3-$\\
$956\,527\,192\,983\,879\,278\,749\,912\,919\,984\,460\,784\,277\,308\,422^{\scriptscriptstyle{\wedge}}3 =
8\,606\,991\,384\,578$

\noindent So 
\begin{align*}
& n=\text{IntegerPart}\left[ \frac{2867531822072415582339770335128768243837513448573364841260551}{8606991384578} \right]= \\ 
& =\,333163087302545313885131123270933026251407963463 \\ 
\end{align*}
Taking in account results obtained we pose the following hypotheses:
\vskip3pt

\noindent {\it Hypothesis 1.}  If $\deg (Q)=0$, then $Q(y)=d^3$ or $Q(y)=2d^3$ ($d$ is a constant). 
\vskip3pt

\noindent {\it Hypothesis 2.}  If $\deg (Q)\neq 0$, then $\sup m=+\infty $.

\bibliographystyle{elsarticle-num}

\end{document}